\documentclass{amsart}
\author{Isabel G. Dotti \and Anna Fino}
\thanks{I.G.Dotti. Research partially supported by
Conicet, SecytU.N.C, Argentina.}
\thanks{A. Fino.  Research partially supported by MURST and GNSAGA (Indam)
of Italy}
\address{
Isabel G. Dotti: FaMAF-CIEM, U. N. de C\'ordoba,\\5000 C\'ordoba,
Argentina\\
Anna Fino: Dipartimento di Matematica, Universit\'a di Torino\\
Via Carlo Alberto 10, 10123 Torino, Italy}
\smallskip

\email{idotti@mate.uncor.edu, fino@dm.unito.it}
\title {Hyperk\"ahler torsion structures invariant by nilpotent Lie groups}
\newtheorem{teo}{Theorem}[section]
\newtheorem{prop}{Proposition}[section]

\newtheorem{lemm}{Lemma}[section]

\newcommand{\beq}{\begin{equation}}
\newcommand{\eeq}{\end{equation}}
\newcommand{\bqn}{\begin{eqnarray}}
\newcommand{\eqn}{\end{eqnarray}}
\newcommand{\bqne}{\begin{eqnarray*}}
\newcommand{\eqne}{\end{eqnarray*}}

\newcommand{\cz}{{\frak z}}
\newcommand{\cv}{{\frak v}}

\newcommand{\cn}{{\frak n}}

\newcommand{\R}{{\Bbb R}}
\renewcommand{\H}{{\Bbb H}}

\begin{document}
\maketitle
\begin{abstract}
We study HKT structures on nilpotent Lie groups and on associated
nilmanifolds. We exhibit three weak HKT structures on $\R^8$ which are
homogeneous with respect to extensions of Heisenberg type Lie groups. The
corresponding hypercomplex structures  are of a special kind, called
abelian.  We prove that on any 2-step nilpotent Lie group all invariant HKT
structures arise from abelian hypercomplex structures.  Furthermore, we use
a correspondence between abelian hypercomplex structures and subspaces of
${\frak sp}(n)$ to produce continuous families of compact and noncompact of
manifolds carrying non isometric HKT structures. Finally, geometrical
properties of invariant HKT structures on 2-step nilpotent Lie groups are
obtained.

\end{abstract}
\section{Introduction}

 Metric connections having totally skew-symmetric torsion arise in a
natural way in
theoretical and mathematical physics. For example,
the geometry of such connections is present on the target space of
supersymmetric sigma models with the Wess-Zumino term
\cite{GHR,Pap,HP2} and , in the supergravity theories, on the moduli space
of a class of black holes \cite{GPS}. Moreover,
the geometry of NS-5 brane solution of type II supergravity theories is
generated by such connection \cite{PT1,PT2,Pa2}.

On any
hermitian manifold $(M,J,g)$ there
exists a unique connection $\nabla$ satisfying $\nabla g = 0, \; \nabla
J = 0$ and
whose torsion tensor
$
c (X,Y,Z) = g (X, T(Y,Z))
$
is
totally skew-symmetric (i.e a three form).
The torsion tensor of this connection is given by $c=- J d J F$, where
$F=g(J.,.)$ is the K\" ahler form for $J$ \cite{Ga}.
  The geometry
of such a connection is called by physicists a KT connection; among
mathematicians this
connection is known as the Bismut connection \cite{Bi}.

 Let $M$ be a smooth manifold  with a hypercomplex
structure $\{ J_i\}_{i=1,2,3}$ and a
riemannian metric $g$.  $M$ is said to be a hyperhermitian manifold if it
is hermitian
with respect to every
$J_i$, $i = 1,2,3$.  A given hyperhermitian manifold $(M, \{ J_i\}_{i=1,2,3},
 g)$  is an HKT (hyperk\"ahler torsion) manifold (\cite{Pap}) if there
 is a
connection $\nabla$ such that
\begin{equation}\label{eq:0}
\nabla g = 0, \quad \nabla J_i = 0, i = 1,2,3, \quad  c (X,Y,Z) = g (X,
T(Y,Z))\;\mbox{a three form.}
\end{equation}
 Such a connection is known as
 an
HKT connection
 in physics literature; its geometry
is known as an HKT geometry.  HKT structures are called strong or weak
depending on whether the torsion $c$ is closed or not. Due to the
uniqueness of the Bismut
connection, a hyperhermitian manifold $M$
will admit an HKT connection  if and only if
$J_1 d J_1 F_1 = J_2 d J_2 F_2 = J_3 d J_3 F_3$ (where $F_i,\; i=1,2,3$ is
the K\"
ahler form associated to  $(J_i, g)$) or equivalently if
$\overline{\partial}_{J_1}(F_2 - i F_3)=0$ \cite{GP}. If this connection
 exists, it is unique \cite{Ga}. Moreover, by \cite{IP}  the associated
  Lee forms $\theta_i=
J_i d^* F_i$ coincide for $i=1,2,3.$

Every $4$-dimensional hyperhermitian manifold is HKT.
If the dimension is $8$ we obtained in \cite{DF1} all simply connected
nilpotent Lie groups which carry invariant abelian
hypercomplex structures.  There are three such groups and they are central
extensions of Heisenberg type Lie groups (see
Example 1 in Section 2).  The abelian hypercomplex structures give rise to
weak HKT structures on these groups (see
Proposition 2.1), with respect to any compatible and invariant riemannian
metric.  These groups are diffeomorphic to $\R^8$.
In coordinates $(x_1,...,x_4, y_1,...,y_4)$ the corresponding HKT metrics
are given by:
    $$\begin{array}{cll}
 g_1&=&\sum dx_i^2 + (dy_1-\frac{1}{2}(x_1dx_2-x_2dx_1-x_3dx_4+x_4dx_3))^2
+\sum_{j\geq 2} dy_j^2,\\  g_2&=&\sum dx_i^2 + dy_1^2 +(dy_2-\frac{1}{2}
(x_1dx_3-x_3dx_1+x_2dx_4-x_4dx_2))^2+\\& & (dy_3-\frac{1}{2}
(x_1dx_4-x_4dx_1-x_2dx_3+x_3dx_2))^2+ dy_4^2,\\
g_3&=&\sum dx_i^2 + (dy_1-\frac{1}{2}(x_1dx_2-x_2dx_1-x_3dx_4+x_4dx_3))^2+
 \\ & & (dy_2-\frac{1}{2}(x_1dx_3-x_3dx_1+x_2dx_4-x_4dx_2))^2+
\\& &(dy_3-\frac{1}{2}(x_1dx_4-x_4dx_1-x_2dx_3+x_3dx_2))^2+
dy_4^2.\end{array}$$
These metrics have a transitive nilpotent group of isometries (hence they
are complete) and they are non isometric to each other.

 The 8-dimensional HKT structures  obtained above are
 associated to abelian hypercomplex structures.
One of the main goals of this paper is to prove that on any   2-step
nilpotent Lie groups all invariant HKT structures arise this way (see
Theorem 3.1).

On the other hand, the correspondence given in \cite{ba} between abelian
 hypercomplex structures on 2-step nilpotent Lie groups and
subspaces of ${\frak sp}(n)$, gives a method to construct infinitely many
compact and non compact families of manifolds carrying non isometric HKT
structures.
By using
 this construction we
 show in Section 4 that there exist  non trivial deformations of
homogeneous HKT structures on ${\R}^{4l},\; l\geq 3$.  Moreover,
for rational parameters one obtains infinitely many HKT compact quotients
of nilpotent Lie groups  by discrete subgroups.  This is in contrast with
results in \cite{BG},\cite{H} in the K\"ahler case.

In the last section we analyze some geometrical properties of invariant
 HKT structures on 2-step nilpotent Lie groups.  We show that in this class,
and with respect to the Bismut connection, the Ricci tensor
  is symmetric, hence by \cite{IP} the torsion 3-form $c$
 is co-closed, every one form in the dual of the center is
 parallel and all
 Lee forms are zero.  This last assertion says that the corresponding
 riemannian manifolds are hermitian semik\" ahler \cite{GH} or hermitian
balanced  \cite{GI}.  These seem to be the first examples of this type.
In the particular case of dimension $8$, using the  explicit
 description of the Bismut connection we show that its Ricci tensor has two
distinct
 eigenvalues $(0, -\lambda, \lambda >0)$ and only one of the groups
carrying invariant HKT structure has parallel torsion.

The authors wish to thank S. Ivanov for useful conversations on the subject
of this paper.

\section{Hyperk\"ahler torsion structures on groups}

A hypercomplex structure on a Lie algebra $\mathfrak g$ is a triple of
 endomorphisms $\{J_i\}_{i=1,2,3}$ satisfying the quaternion relations
$J_i^2= -I,\; i=1,2,3,\;\; J_1J_2=-J_2J_1=J_3$,
together with the vanishing of
 the Nijenhuis tensor
$$N_i(X,Y)= J_i([X,Y]-[J_iX,J_iY])-([J_iX,Y]+[X,J_iY]),$$
where $X,Y\in \frak g$ and $i=1,2,3$.

The hypercomplex structure will be called abelian if
$$[J_i X, J_i Y]= [X,Y],$$
for all $X, Y \in \frak g, \;\; i=1,2,3. $
Abelian hypercomplex structures were previously considered in \cite{BD},
 \cite{DF1}, \cite{DF2}; they can only occur on solvable Lie algebras
 (\cite{bilbao}).

Let $\frak g$ be a Lie algebra endowed with a hypercomplex
structure $\{J_i\}_{i=1,2,3}$ and an inner product $g$,
compatible with the hypercomplex structure, that is
\begin{eqnarray}\label{eq:2}
g (X,Y) = g (J_1X,J_1Y) = g (J_2X,J_2Y)= g(J_3X,J_3Y) .
\end{eqnarray}
for all $X,Y \in \frak g.$   Assume furthermore that the hypercomplex
structure together with the inner product given on
 $\frak g$ satisfy the extra condition
\begin{eqnarray}\label{eq:3} g([J_1X,J_1Y],Z)+g([J_1Y,J_1Z],X)+g([J_1Z,J_1X],Y
)=
\nonumber\\ g([J_2X,J_2Y],Z)+g([J_2Y,J_2Z],X)+
g([J_2Z,J_2X],Y)=\\
g([J_3X,J_3Y],Z)+g([J_3Y,J_3Z],X)+g([J_3Z,J_3X],Y),\nonumber
\end{eqnarray}
for all $X,Y,Z \in \frak g$.
Note that if  one substitutes $X,Y,Z$ by $J_3X, J_3Y, J_3Z$ in the second
and third row of (\ref{eq:3}), then one obtains
\begin{eqnarray}\label{eq:4}
g(J_3[X,Y],Z)+g(J_3[Y,Z],X)+g(J_3[Z,X],Y)=\nonumber\\
g(J_3[J_1X,J_1Y],Z)+g(J_3[J_1Y,J_1Z],X)+g(J_3[J_1Z,J_1X],Y),
\end{eqnarray}
and conversely, (\ref{eq:4}) implies that the last two rows in
(\ref{eq:3}) are equal.  Since
the equality of any two rows
in (\ref{eq:3}) gives equality of the three rows, one has  in particular that
 (\ref{eq:3}) and
(\ref{eq:4}) are equivalent.

An HKT structure $(\{J_i\}_{i=1,2,3},\;g)$ on a Lie algebra ${\frak g}$
consists of a hypercomplex structure  $\{J_i\}_{i=1,2,3}$  toghether
with an inner product $g$ on ${\frak g}$
satisfying conditions (\ref{eq:2}) and (\ref{eq:3}) or conditions
(\ref{eq:2}) and (\ref{eq:4}).
If $G$ is a Lie group with Lie algebra $\frak g$ carrying an HKT structure, by
left translating the $J_i,\; i=1,2,3$ and the inner product $g$,
one obtains in $G$  an invariant HKT structure. Indeed, in
this case one finds that the Bismut connection is defined by the equation
\begin{eqnarray}\label{eq:5}
g(\nabla_XY,Z)=\tfrac{1}{2}\{g([X,Y]-[J_iX,J_iY],Z)\\
-g([Y,Z]+[J_iY,J_iZ],X)+g([Z,X]-[J_iZ,J_iX],Y)\}\nonumber.
\end{eqnarray}
for $X,Y,Z$ left invariant vector fields.  A verification shows that this
connection satisfies  (\ref{eq:0}).

When the hypercomplex structure is abelian, (\ref{eq:3})
 is always satisfied and moreover, the HKT structure is weak.  Indeed,
to prove the last assertion, we note that the Bismut connection $\nabla$ and
its torsion $c$ are given by
$$g(\nabla_XY,Z)= -g([Y,Z],X),$$
$$c(X,Y,Z)= g(X,T(Y,Z))=
(-1)(g([X,Y],Z)+g([Y,Z],X)+g([Z,X],Y)).$$
  Since
\begin{eqnarray}
dc(X,Y,Z,W)= -2g([X,Y],[Z,W])\nonumber\
 +2g(([X,Z],[Y,W])- 2g([X,W],[Y,Z])\nonumber,\end{eqnarray}
one obtains in particular
\begin{eqnarray}\label{eq:6}
dc(X,J_1X,J_2X,J_3X)= 2||[X,J_1X]||^2+ 2||[X,J_2X]||^2+ 2||[X,J_3X]||^2
\end{eqnarray}
and
\begin{eqnarray}\label{eq:7}
dc(X,J_1X,Y,J_1Y)= -2g([X,J_1X],[Y,J_1Y])\nonumber\\
 + 2g([X,Y],[X,Y])+ 2g([X,J_1Y],[X,J_1Y]).
\end{eqnarray}
Equations (6) and (7) imply that $dc\neq 0$ unless $\frak g$ is abelian.
Indeed, the condition $dc=0$
in (6) implies $g([X,J_1 X], [X, J_1 X])=0$ and
substituting this in (7) gives $g( [X,Y], [X,Y] )=0$.
 Hence, we have proved

\begin{prop} Every abelian hypercomplex structure on a non abelian Lie
group $G$ give rise to an invariant weak HKT structure on $G$.
\end{prop}

\subsection {Examples}

\begin{enumerate}
\item
Let $H_i(n)$ for $1\le i \le 3$, denote respectively the real,
 complex or quaternionic Heisenberg groups. The  hypercomplex
 structures on the 8-dimensional nilpotent Lie groups  $N_1={\R}^3\times
 H_1(2),\; N_2={\R}^2\times H_2(1)\; N_3={\R}^1\times H_3(1)$ constructed in
 \cite{DF1} are abelian
and give rise, with respect to any compatible and invariant riemannian
metric, to weak HKT-structures on these groups.  Moreover, as proved in
\cite{bilbao} the
Obata connection associated to any hypercomplex structure on $N_i$, $i =
1,2,3$ is flat.
 Their corresponding Lie algebras are  $\cn_i =
 \cv \oplus \cz, \; i=1,2,3$,  $\cv = \mbox {span} \{ e_1,
e_ 2, e_3, e_4\}$ and $\cz = \mbox {span} \{ e_5, e_6,
e_7, e_8 \}$, with non zero brackets
$$\lbrack e_1, e_2 \rbrack =
-\lbrack e_3, e_4 \rbrack =  e_5$$ in $\cn_1$,
$$\lbrack e_1, e_3 \rbrack =  \lbrack e_2, e_4 \rbrack = e_6\,;\;\;
\lbrack e_1, e_4 \rbrack =- \lbrack e_2, e_3 \rbrack =
e_7$$ in $\cn_2$ and
$$\lbrack e_1, e_2 \rbrack = - \lbrack e_3,
e_4
\rbrack = e_5; \quad\lbrack e_1, e_3 \rbrack = \lbrack e_2, e_4 \rbrack
= e_6; \quad \lbrack e_1, e_4 \rbrack = - \lbrack e_2, e_3 \rbrack =
e_7$$ in $\cn_3.$  The hypercomplex structure is given by $J_ie_1=e_{i+1},
 J_ie_5=e_{5+i},\; i=1,2,3$, $J_i^2=-I,\; J_1J_2=-J_2J_1=J_3$.
We note that
 these nilpotent Lie groups do admit lattices,
hence we also obtain compact examples.

\item The 8-dimensional Lie group $H_1(2)\times SU(2)$ has an invariant
 weak HKT structure (see \cite{poon}). This group does not admit invariant
 abelian hypercomplex structures since it is not solvable (\cite{bilbao})

\item The  12-dimensional  3-step nilpotent Lie group with non zero brackets
$$
\begin{array}{l}
\lbrack e_1, e_2 \rbrack = - [e_5, e_6] = - e_{10}, [e_2, e_5] = - [e_1,
e_6] = - e_{11},\\
\lbrack e_1, e_4 \rbrack =  [e_2, e_{10}] = [e_5, e_8] = [e_6, e_{11}] = -
e_{12}
\end{array}
$$
admits an abelian hypercomplex structure  $\{ J_i \}_{i = 1,2,3}$ given by
(\cite{bilbao})
$$
\begin{array} {l}
J_1 e_1 = e_2, J_1 e_3 = e_{12}, J_1 e_4 = e_{10}, J_1 e_5 = e_6, J_1 e_7 =
e_9, J_1 e_8 = e_{11},\\
J_2 e_1 = e_6, J_2 e_2 = e_{5}, J_2 e_3 = e_{9}, J_2 e_4 = e_{11}, J_2 e_8
= -e_{10}, J_2 e_7 = -e_{12}
\end{array}
$$
whose  associated Obata connection is
not flat, by \cite{bilbao}. On the other hand
 the hypercomplex structure together with the metric such that  the
above basis is orthonormal give a weak invariant HKT
structure.
 \end{enumerate}

\section{HKT structures on nilpotent Lie groups}

In this section we will restrict to the case of invariant HKT  structures  on
 nilpotent Lie groups.

Let $\mathfrak n$ be a nilpotent
 Lie algebra, that is a Lie algebra satisfying
${\mathfrak n}^k = 0$ for some $k \geq 1$, where ${\mathfrak n}^i$ is the
chain of
ideals defined  inductively by ${\mathfrak n}^0 = \mathfrak n$
and
$$
{\mathfrak n}^i = \lbrack {\mathfrak n}^{i - 1}, {\mathfrak n} \rbrack,
\quad i \geq 1.
$$
One says that $\mathfrak n$ is $k$-step nilpotent if ${\mathfrak n}^k = 0$
 and
${\mathfrak n}^{k - 1} \neq
0$, $k \geq 1$.

Let $\{J_i\}_{i=1,2,3}$ be a hypercomplex structure on a nilpotent Lie
 algebra $\frak n$.  In order to prove the main result of this section we
first prove two useful
 lemmas.

\begin{lemm} If $\{J_i\}_{i=1,2,3}$ is a hypercomplex structure on an
 $s$-step nilpotent Lie algebra $\frak n$ then the
inclusion

\begin{equation}\label{eq:8}
{\mathfrak n}^{s-1}_Q =
{\mathfrak n}^{s-1} + J_1 {\mathfrak n}^{s-1} + J_2 {\mathfrak n}^{s-1} +
 J_3
{\mathfrak n}^{s-1} \subset
{\mathfrak n}
\end{equation}
is proper.
\end{lemm}

\noindent {\bf Proof.}  Suppose it is not. Let $X \in {\mathfrak n}$
and write $X =
X_0 + J_1 X_1
+ J_2 X_2 + J_3 X_3$, $X_i \in {\mathfrak n}^{s-1}, i=0,1,2,3 $. If $Y
\in \mathfrak n$, then
$$
\lbrack X, Y \rbrack =  \lbrack J_1 X_1, Y \rbrack +
\lbrack J_2 X_2, Y \rbrack + \lbrack J_3 X_3, Y \rbrack
,$$
since $\mathfrak n$ is $s$-step. Write $Y=
Y_0 + J_1 Y_1
+ J_2 Y_2 + J_3 Y_3$, $Y_i \in {\mathfrak n}^{s-1}, i=0,1,2,3 $
and substitute in the previous expression obtaining
\begin{eqnarray}\label{eq:9}
[X,Y]& = &\sum_{i=1}^{i=3}[J_iX_i, J_iY_i] + [J_1X_1, J_2Y_2+J_3Y_3]+
\nonumber \\
& & [J_2X_2, J_3Y_3+J_1Y_1]
+ [J_3X_3, J_1Y_1+J_2Y_2].
\end{eqnarray}
 Denote by $\frak z$ the center of $\frak n$ (note that ${\mathfrak n}^{
s-1}\subset \frak z$) and observe that the
integrability of
$J_l$ gives, for $l = 1,2,3$,
\begin{eqnarray}\label{eq:10}
[J_l U, J_l V] =0,\;
 \;\; {\rm  if}\;\; U,V \in {\mathfrak z},\\
\label{eq:11} J_l[J_l U,  V]= [J_l U, J_l V] \;{\rm if}\;\; U \in {\frak
z}\;\;{\rm
and}\;\; V \in {\mathfrak n}.
\end{eqnarray}
Note that (11) above implies
\begin{eqnarray}\label{eq:12}[J_1U,J_3V]=J_1[J_1U,J_2V],\;\;
[J_2U,J_3V]=J_2[J_2U,J_1(-V)]
, \;\; U,V \in {\mathfrak z}.
\end{eqnarray}
Hence
using (10) and (12)  in the expression of $[X,Y]$ and setting
${\mathfrak u}= [J_1 {\mathfrak n}^{s-1}, J_2 {\mathfrak n}^{s-1}]$
one obtains
$$ {\mathfrak n}^1 = {\mathfrak u} + J_1 {\mathfrak u} +J_2{\mathfrak u}.$$
Then
$${\mathfrak n} = {\mathfrak n}^{s-1}_Q \subset
{\mathfrak n}^{1} + J_1 {\mathfrak n}^{1} + J_2 {\mathfrak n}^{1} + J_3
{\mathfrak n}^{1} \subset
{\mathfrak u} + J_1 {\mathfrak u} + J_2 {\mathfrak u}
+J_3{\mathfrak u}$$
and as a consequence ${\mathfrak n}= {\mathfrak n}^1 + J_3{\mathfrak n}^1$
 contradicting the fact that for any
invariant complex structure on a nilpotent Lie group there exists a closed
$(1,0)$ form (\cite{sal}).

\medskip

\noindent{\bf Remark 1.}
We observe that ${\mathfrak z}_Q =
{\mathfrak z} + J_1 {\mathfrak z} + J_2 {\mathfrak z} + J_3
{\mathfrak z}$
 can be all of $\frak n$.  Indeed in
\cite{DF2} such an example is given of a  hypercomplex nilpotent Lie
algebra of dimension 8 having a 5-dimensional center.
Also, when the hypercomplex structure is abelian one
 has (see \cite{DF2}) that the inclusion
$$
{\mathfrak n}^i_Q =
{\mathfrak n}^i + J_1 {\mathfrak n}^i + J_2 {\mathfrak n}^i + J_3
{\mathfrak n}^i \subset
{\mathfrak n}^{i - 1} + J_1 {\mathfrak n}^{i - 1} + J_2
 {\mathfrak n}^{i - 1}
+ J_3 {\mathfrak n}^{i - 1} = {\mathfrak n}^{i - 1}_Q ,$$ $i\geq 1$,
is proper.

\begin{lemm}
Let $\frak n$ be a 2-step nilpotent Lie algebra with an HKT structure.
 If ${\mathfrak n}^{1}_Q= {\mathfrak n}^{1} + J_1 {\mathfrak n}^{1} +
 J_2 {\mathfrak n}^{1} + J_3
{\mathfrak n}^{1}$ and
  $\frak m =({\mathfrak n}^{1}_Q)^{\perp}$ then

i) $[{\mathfrak n}^{1}_Q,{\mathfrak n}^{1}_Q]=0.$

ii) $[{\mathfrak n}^{1}_Q,{\mathfrak m}]=0.$

\end{lemm}
\noindent{\bf Proof.} Note first that by the previous lemma the subspace
 ${\mathfrak n}^{1}_Q$ is proper.  To prove assertion i) we observe
 first that ${\mathfrak n}^{1}
+ J_l {\mathfrak n}^{1}$ is an abelian subalgebra, for any $l=1,2,3$,
 since
${\mathfrak n}^{1}$ is contained in the center of ${\frak n}$ (see (10)).
 We
show next that $[J_1{\mathfrak n}^{1}, J_2{\mathfrak n}^{1}]=0.$  Take
$X\in {\mathfrak n}^{1},\; Y=J_3Y',\;Z=J_3Z',\; Y',Z'\in {\mathfrak n}^{1}$
 and substitute
in identity (\ref{eq:4}).  Then one obtains
$$0=g(J_3[J_1X,-J_2Y'],J_3Z')+ g(J_3[-J_2Z',J_1X],J_3Y') $$
or equivalently $ad_{J_1X}J_2$ is a self-adjoint transformation of
 ${\frak
n}^{1}$.  Since its square is
zero by  (11), $[J_1X,J_2Y]=0, X,Y\in {\mathfrak n}^{1}.$  Using (12) one
shows $[J_iX,J_jY]=0, X,Y\in {\mathfrak n}^{1}$
in the remaining cases.
\medskip
We next prove  assertion ii).
Substituting $X,\;Z \in {\mathfrak n}^{1}$ and $Y \in {\mathfrak m}$ into
the last two lines
 of (\ref{eq:3}) implies that $ad_{J_2 Y} J_2 -ad_{J_3 Y} J_3$ is a
self-adjoint transformation of ${\mathfrak n}^{1}$.
 Since it is also nilpotent (see (11)) it must be zero hence
\begin{eqnarray}
[J_2X,J_2Y]=[J_3X,J_3Y],\;\;X \in {\mathfrak n}^{1},\;\; Y\in {\mathfrak m}.
\end{eqnarray}
Take now $X, \tilde{Z} \in {\mathfrak n}^{1},\;\; Y\in {\mathfrak m},\;
Z=J_3\tilde{Z}$ and substitute in (\ref{eq:4}) obtaining
\begin{eqnarray}
g([J_1X,J_1Y],\tilde{Z})= g(J_3([Y,J_3\tilde{Z}] + [J_1Y,J_2\tilde{Z}]),X).
\end{eqnarray}
Since $[J_2\tilde{Z}, J_1Y]= [J_2\tilde{Z}, J_2(J_1J_2Y)]=
[J_3\tilde{Z}, J_3(J_1J_2Y)]= - [J_3\tilde{Z}, Y]$ by (\ref{eq:11})
 it follows that
$[J_1X,J_1Y]=0$ or equivalently
 $[J_1{\mathfrak n}^{1}, {\mathfrak m}]=0$.  The proof in the remaining
cases is similar.

\begin{teo}
 The hypercomplex structure of an invariant HKT structure
 on any 2-step nilpotent Lie group is  abelian.
\end{teo}

\noindent{\bf Proof.}
From i) and ii) of the  previous lemma  it follows
that ${\mathfrak n}^1_Q \subset {\frak z}$
where $\frak z$ stands for the center of $\frak n$. If $Y,Z\in \frak m$ then
(\ref{eq:3}) implies $[J_1Y,J_1Z]=[J_2Y,J_2Z]=[J_3Y,J_3Z]$.  Since
$\frak m$ is
$J_i-$invariant, it follows that for
any $Y,Z \in \frak m$ then $[Y,Z]=[J_iY,J_iZ],\; i=1,2,3.$
Now
it is straightforward to show that the hypercomplex structure is abelian, by
decomposing any given $U,V \in \frak n$,
as $U=U_1+U_2,\;V=V_1+V_2 $ according to $\frak n =
{\mathfrak n}^1_Q \oplus \frak m .$

\bigskip

\noindent{\bf Remark 2.} If we restrict to the 8-dimensional case  it was
proved in \cite{DF1} that the only 2-step nilpotent Lie groups carrying
abelian hypercomplex structures are the groups $N_i, i=1,2,3$ considered in
Example 1.
As remarked in \cite{DF1} $N_1, N_2$ can only carry abelian hypercomplex
structures. On the other hand, $N_3$ does admit
non-abelian hypercomplex structures, thus by Theorem 3.1,  $N_3$
endowed with a non abelian hypercomplex
structure admits no invariant metric such that it becomes an HKT manifold.

\medskip

\section{Deformations of HKT structures}

  In \cite{ba} M. L. Barberis proved that there is a one to one
correspondence between  injective linear maps
$j:{\R}^m\rightarrow {\frak sp(k)}\;(m\leq k(2k+1))$ and
 2-step nilpotent Lie
algebras $\frak n$ with
dimension  $([\frak n, \frak n])=k$ carrying an abelian hypercomplex
structure.   Using Theorem 3.1
 one can rephrase the above result saying that the correspondence is
between injective linear maps $j:{\R}^m\rightarrow {\frak
sp(k)}\;(m\leq k(2k+1))$ and  2-step nilpotent Lie algebras $\frak n$ with
dimension $([\frak n, \frak n]=k)$ carrying an HKT  structure.

 We reproduce the construction given in \cite{ba} with the only
modification introduced by Theorem 3.1.

Let $\frak n$ be a 2-step nilpotent Lie algebra  with an HKT structure
$(\{J_i\}_{i=1,2,3},g)$ and consider the orthogonal
decomposition $\frak n = \frak v \oplus \frak z$, with $\frak z$ the center
of $\frak n$.  Note that since $\frak n$ is 2-step
nilpotent one has that $[\frak n,\frak n]\subset \frak z$.
 Since $\{J_i\}_{i=1,2,3}$ is an abelian hypercomplex structure then it
preserves $\frak z$, hence $\frak v$, since it is hyperhermitian with
respect to $g$.  Let $j: [\frak n,\frak n]\rightarrow
{\frak so(v)}$ be defined by
$$ g(j_zX,Y)=g([X,Y],z),\;\; X,Y \in {\frak v},\; z\in \frak z.$$
Then $j$ is  one to one and if ${\frak sp(v)}= \{T\in {\frak so(v)}:
TJ_i=J_iT, i=1,2,3\}$ one has $j_z \in {\frak sp(v)}$
for all $z\in [\frak n,\frak n]$ since
$$ g(j_zJ_iX,J_iY)=g([J_iX,J_iY],z)=g([X,Y],z)=g(j_zX,Y).$$
Conversely, given $j: {\R}^m\rightarrow {\frak sp(k)}$, fix $0\leq s\leq 3$
with $s+m
\equiv 0$ mod($4$) and set $\frak n ={\R}^k\oplus {\R}^s\oplus {\R}^m$ with
 $g$ the
canonical inner product. Define the bracket such
that ${\R}^s\oplus {\R}^m$ is central and  $g([X,Y],z)= g(j_zX,Y)$
 if $X,Y \in
{\R}^k, \;z\in {\R}^m$.  Let $\{J_i\}_{i=1,2,3}$ be the
endomorphisms of ${\R}^k$ defining
${\frak sp(k)}$ extended to all of $\frak n$ by anticommuting complex
endomorphisms on ${\R}^s\oplus {\R}^m$ compatible
with the metric.  It is easy to verify that the resulting hypercomplex
structure is abelian, hence it is an HKT structure on
$\frak n$.

\medskip

\noindent{\bf Remark 3.} By applying the previous construction to the case
$m=1$ and   $j_z$
any complex structure commuting with
the complex structures defining ${\frak sp(k)}$ (for a fixed $z\neq 0$ in
$\R$), the resulting
algebra is an extension of the Heisenberg algebra and the HKT-structure
is that obtained in \cite{GP}(5.2).

\medskip

We next give some non trivial
deformations of HKT structures.

\medskip

Fix in ${\R}^{4l}, l\geq 2$, identified with ${\H}^l$, $\H$ the quaternions,
the hypercomplex
structure $\{J_1,J_2,J_3 \}$
given by right multiplication by $(i,...,i),\; (j,...,j),$
and $(-k,...,-k)$ respectively.
Let $t > 0$ and $j^t : {\R}^2 \rightarrow
{\frak sp(l)}$, with
$${\frak sp(l)}= \{T\in {\frak so(4l)}: TJ_i=J_iT,
i=1,2,3\},
$$
 be given by
$$j^t_{e_1}= L_{(i,...,i,i)}, \;\; j^t_{e_2}= L_{(j,..,j,tj)}$$ where $L$
stands for
left multiplication, and $e_1,e_2$
denotes a basis of ${\R}^2$.  It is clear that $j^t$ is a mapping into
 ${\frak
{sp}}(l)$ since left and right multiplication
commute.

Similarly, if $l\geq 3$, let $(t,s)$ be such that
$0<t<s<1 $ and $j^{t,s} : {\R}^3 \rightarrow
{\frak sp(l)}$, be given by
$$j^{t,s}_{e_1}= L_{(i,...,i,i)}, \;\; j^{t,s}_{e_2}= L_{(j,...,tj,j)},
\;\; j^{t,s}_{e_3}=
L_{(k,...,k,sk)},$$ where $L$ stands for left multiplication,
and $e_1,e_2, e_3$ denotes a basis of ${\R}^3$.

Let ${\frak n}_t $(resp. ${\frak n}_{t,s}$) = ${\R}^{4l}\oplus {\R}^4$ be
the 2-step nilpotent
 Lie
algebra with the HKT structure constructed as above and
let $N_t$ (resp. $N_{t,s}$) be the simply connected Lie group with Lie
algebra ${\frak n}_t$ (resp. ${\frak n}_{t,s}$)
and invariant HKT structure induced by left
translating the inner product and hypercomplex structure on ${\frak n}_t$ (resp.
  ${\frak n}_{t,s}$).

\medskip

\noindent{\bf Claim 1. } The riemannian manifolds $N_t$ and $N_{t'}$ (resp.
$N_{t,s}$ and $N_{t',s'}$) are
isometric
if and only if $t=t'$ (resp. $(t,s)=(t',s')$).

\medskip

According to E. Wilson \cite{W}, if two nilpotent Lie groups with left
invariant metric are isometric
 there exists an isomorphism which is also an isometry.  Hence, its
derivative is an orthogonal Lie algebra isomorphism
between the corresponding Lie algebras. Assume then that $f$ denotes an
orthogonal isomorphism from ${\frak n}_t$ onto
${\frak n}_{t'}.$  Using the description of the Lie brackets given by $j^t$
and $j^{t'}$ respectively, it follows that
${j^t}_{f^{-1}z}= f^{-1}j^{t'}_{z}f$ for all $z \in {\R}^2.$  Squaring both
sides in the previous equality and carrying out a
tedious but straightforward computation one finds that $t=t'$. Using
similar arguments one can show that the
riemannian manifolds $N_{t,s}$ and $N_{t',s'}$ are
isometric if and only if $(t,s)=(t',s')$.

\medskip

\noindent{\bf Claim 2. } The riemannian manifolds $N_q, \; q\in {\mathbf
Q},\; q >0$,  do admit  discrete subgroups $\Gamma_q$ such that $T_q=
N_q/\Gamma_q$ is compact.  Furthermore, the (non homogeneous) $T_q$  with
the induced metrics, are not isometric to eacch other for differents $q$'s.

\medskip

We recall that according to \cite{ma} a nilpotent Lie group $N$ admits a
discrete subgroup $\Gamma$ such that $N/\Gamma$ is compact if and only if
its Lie algebra $\frak n$ admits a basis with rational structure constants.
But this is clear in the case
$t$ rational in the definition of ${\frak n}_t.$   Furthermore, an isometry
between $T_q$ and $T_q'$ lifts to an isometry between $N_q$ and $N_q',$
which is impossible by Claim 1.

\medskip

A compact quotient of a nilpotent Lie group $N$ by a discrete subgroup is
called a {\em nilmanifold}.

\begin{teo}There exists a one parameter (resp. two parameter family) of
homogeneous HKT structures on ${\R}^{4l}, l\geq 3$ ( resp. ${\R}^{4l},
l\geq 4$).  Moreover, for rational parameters
there exists infinitely many non isometric HKT nilmanifolds.
\end{teo}

\noindent{\bf Proof.} To prove
the existence of a one parameter family of non isometric HKT structures on
${\R}^{4l}, l\geq 3$(respectively a two parameter family of HKT structures on
 ${\R}^{4l}, l\geq 4$) one
uses the fact that the exponential map $exp_t$ (resp.
$exp_{t,s}$) is a diffeomorphism from ${\R}^{4l} \rightarrow N_t$ (resp.
${\R}^{4l} \rightarrow N_{t,s}$).
 The pullback by $exp_t$ (resp. $exp_{t,s}$) of the invariant HKT-structures
on $N_t$ (resp. $N_{t,s})$ together with Claim 1 gives the asserted
deformation.   The second statement in Theorem 4.1 follows from Claim 2.

\medskip

\noindent{\bf Remark 4.} The existence of infinitely many nilmanifolds
carrying HKT structures is in contrast with the K\"ahler case (compare
\cite{BG}, \cite{H}).

\section{Geometrical consequences}

Let $\frak n$ be a 2-step nilpotent Lie algebra
with an HKT structure  $(\{J_i\}_{i=1,2,3}, g)$, $\nabla$ the Bismut
connection and $c$ the torsion 3-form. We show next that the Ricci
tensor of $\nabla$  is symmetric (hence $c$ is co-closed by
\cite{IP}) and the Lee forms are zero, hence the corresponding
 riemannian manifolds are hermitian semik\" ahler  (according to
\cite{GH}) or hermitian balanced  \cite{GI}.

 In Section 2 we observed that the Bismut connection
 associated to an abelian
hypercomplex structure and its torsion 3-form
 were given respectively  by
$$ \begin{array}{l}
g(\nabla_XY,Z)= -g([Y,Z],X),\;\;\\
c(X,Y,Z)=
(-1)(g([X,Y],Z)+g([Y,Z],X)+g([Z,X],Y)),
\end{array}$$  for any $X,Y,Z \in \frak n$.

Decompose $\frak n= \frak v \oplus \frak z$ where
$\frak z$ is the center of $\frak n$ and $\frak v$ its orthogonal
 complement.  It follows easily that
\begin{enumerate}
\item[i)]$\nabla_XZ=0,\;\; Z\in \frak z, X\in \frak n.$
\item[ii)]$\nabla_VX=0,\;\; V\in \frak v, X\in \frak n.$
\item[iii)]$\nabla_XY\in \frak v ,\;\; X,Y\in  \frak n.$

\end{enumerate}
The curvature tensor associated to $\nabla$ is
 $R(X,Y)=[\nabla_X,\nabla_Y]-\nabla_{[X,Y]}.$
  As a consequence of i) and ii) above it follows that
\begin{enumerate}
\item[iii)]$g(R(X,Z)Y,Z')=0,\;\; Z, Z'\in \frak z, X, Y\in \frak n.$
\item[iv)]$g(R(X,V)Y,V')= g([X,V],[Y,V']),\;\;  V, V'\in \frak v,
 X,Y\in \frak n.$
\end{enumerate}
In particular, the  Ricci tensor $\rho$  of $\nabla$ is symmetric,
 given by
$$
\rho (X, Y) = \sum_j g (R(X, V_j) Y, V_J) =  \sum_j g([Y,V_j],[X,V_j]),$$
 where $ V_j$ is an  orthonormal basis of $ \frak v$. It is worth
 to point out that the Ricci tensor of the Bismut connection
 $\nabla$ is not symmetric in general. By \cite[Corollary 3.2]{IP}
 the Ricci tensor
of $\nabla$ is symmetric if and only the torsion 3-form $c$ is co-closed.
  In particular,  on a 2-step nilpotent Lie algebra with an
HKT structure, the torsion 3-form $c$ is always co-closed.
Moreover, it follows from iii) that for any 1-form $\alpha$ in
 the dual of the center $\frak z^*$ is parallel with respect to
 the Bismut connection, thus giving a reduction of its holonomy group.

One can also verify, using  the expression of $\rho$ above
that
\begin{enumerate}
\item[v)]$\rho(Z,X)=0,\;\; Z\in \frak z, X\in \frak n.$
\item[vi)]$\rho (V,J_lV)=0, V\in \frak v,\; l=1,2,3.$
\item[vii)]$\rho(V,V)=\rho(J_lV,J_lV), \; V\in \frak v\;, l=1,2,3.$
\end{enumerate}
The last two assertions follow from
$$\rho(V,J_lV)=\sum g([V,V_j],[J_lV, V_j])=
- \sum g([J_lV,J_lV_j],[V, J_lV_j])=-\rho(V,J_lV),$$
 and
$$\rho(V,V)=\sum g([V,V_j],[V, V_j])=
\sum g([V,J_lV_j],[V, J_lV_j])=\rho(J_lV,J_lV).$$

Finally, to show that the Lee forms are trivial we need to recall that by
 \cite{Ga,IP}
$$
\theta (X) = - 1/2 \sum_{i=1}^{2n} c(J_l X,e_i, J_l e_i), \leqno(17)
$$
where $c$ is the torsion 3-form and $\{ e_i \}$ is an orthonormal
 basis of $\mathfrak n$.

In general, one has  that if $Y \in \frak {z}$,
$$
c (J_l X, Y, J_l Y) = 0,
$$
 and if $Y \in \frak v$,
$$
c (J_l X, Y, J_l Y) = - g ([Y, J_l Y], J_l X).
$$
Then, using a basis $V_j, J_1 V_j, J_2 V_j, J_3 V_j$ of
$\frak v$ and letting $l = 1$ in (17)
$$
\theta (X) = - 1/2 \sum_j 2g([V_j,J_1 V_j],J_1 X) + 2g([J_2 V_j,J_3
V_j],J_1 X) =0
$$
since $J_1$ is abelian.

\subsection{Geometry of 8-dimensional examples} We restrict next to
 the case of an 8-dimensional 2-step nilpotent
Lie group.  In \cite{DF1} we showed that the only nilpotent 8-dimensional
 Lie groups
carrying abelian hypercomplex structures were the groups $ N_i, i=1,2,3$
 described in example 2 of 2.1.
The groups $N_i,\; i=1,2,3$ are diffeomorphic to ${\R}^8$ via the inverse
 of the exponential map. Using this diffeomorphism $(x_1,x_2,x_3,x_4,
 y_1,y_2,y_3,y_4): N_i\rightarrow \frak v\oplus \frak z$ as a coordinate
 system
 one can write the three complete HKT metrics on ${\R}^8$ as follows
$$\begin{array}{cll}
 g_1&=&\sum dx_i^2 + (dy_1-\frac{1}{2}(x_1dx_2-x_2dx_1-x_3dx_4+x_4dx_3))^2
+\sum_{j\geq 2} dy_j^2,\\  g_2&=&\sum dx_i^2 + dy_1^2 +(dy_2-\frac{1}{2}
(x_1dx_3-x_3dx_1+x_2dx_4-x_4dx_2))^2+\\& & (dy_3-\frac{1}{2}
(x_1dx_4-x_4dx_1-x_2dx_3+x_3dx_2))^2+ dy_4^2,\\

g_3&=&\sum dx_i^2 + (dy_1-\frac{1}{2}(x_1dx_2-x_2dx_1-x_3dx_4+x_4dx_3))^2+
 \\ & & (dy_2-\frac{1}{2}(x_1dx_3-x_3dx_1+x_2dx_4-x_4dx_2))^2+
\\& &(dy_3-\frac{1}{2}(x_1dx_4-x_4dx_1-x_2dx_3+x_3dx_2))^2+
dy_4^2.\end{array}$$
The metrics are not isometric since they come from non
 isomorphic groups (see \cite{W}).
The  Ricci tensor
associated to the Bismut connection, in these cases is given by
$$\rho(Z,X)=0,\;\; Z\in \frak z, X\in \frak n \quad \rho(V,V)=c, c<0,
\; V\in \frak v, ||V||=1.$$
Indeed, in this case $\dim {\mathfrak v} = 4$ and one can consider
 $\{ V, J_1 V, J_2 V, J_3 V \}$ as basis of $\mathfrak
v$ and apply v),vi),vii) above.

We show next that the torsion 3-form $c$ is parallel with respect to the
 Bismut connection $\nabla$ only in the case of $N_1$.
The Bismut connection $\nabla$  is given on $N_1$  by
$$
\nabla_{e_5} e_1 = - e_2, \nabla_{e_5} e_2 =  e_1,
\nabla_{e_5} e_3 = e_4, \nabla_{e_5} e_4 = - e_3,
$$
on $N_2$ by
$$
\begin{array} {l}
\nabla_{e_6} e_1 = - e_3, \nabla_{e_6} e_2 =  -e_4,
 \nabla_{e_6} e_3 = e_1, \nabla_{e_6} e_4 =  e_2,\\
\nabla_{e_7} e_1 = - e_4,
 \nabla_{e_7} e_2 = e_3,
 \nabla_{e_7} e_3 = - e_2, \nabla_{e_7} e_4 = e_1
\end{array}
$$
and on $N_3$ by
$$
\begin{array} {l}
\nabla_{e_5} e_1 = -e_2, \nabla_{e_5} e_2 =  e_1,
\nabla_{e_5} e_3 = e_4, \nabla_{e_5} e_4 = -e_3,\\
\nabla_{e_6} e_1 = - e_3, \nabla_{e_6} e_2 =  - e_4,
\nabla_{e_7} e_3 = -e_1, \nabla_{e_6} e_4 = e_2,\\
\nabla_{e_7} e_1 = - e_4, \nabla_{e_7} e_2 =  e_3,
\nabla_{e_7} e_3 = - e_2, \nabla_{e_7} e_4 =  e_1,
\end{array}
$$

respectively.

On $N_1$ the torsion 3-form $c$ is given by
$$
e^3 \wedge e^4 \wedge e^5 -  e^1 \wedge e^2 \wedge e^5
$$
and it is parallel with respect to the  Bismut connection.
On $N_2$ the torsion 3-form  $c$ is given by
$$
- e^2 \wedge e^4 \wedge e^6 - e^{1} \wedge e^3 \wedge e^6 + e^{7}
 \wedge e^2 \wedge e^3 - e^{7} \wedge e^1 \wedge e^4
$$
and it is not parallel, since for example $(\nabla_{e_6} c)
(e_1, e_2,
e_7) \neq 0$.
On $N_3$ the torsion 3-form $c$ is given by
$$
 e^{5} \wedge e^2 \wedge e^1 - e^{5} \wedge e^4 \wedge e^3 - e^{2}
 \wedge e^4 \wedge e^6 - e^{1} \wedge e^3  \wedge e^6 +
e^{7} \wedge e^2 \wedge e^3 - e^{7} \wedge e^1 \wedge e^4
$$
and it is not parallel, since for example $(\nabla_{e_6} c) (e_1, e_3,
e_8) \neq 0$.

\bigskip

\noindent
{\bf Concluding remarks}

\smallskip
\noindent
In \cite{GPS} it is shown that the geometry of the moduli space of a class
of black holes in five dimension is HKT and the
relation between the number of supersymmetries of a sigma model and the
geometry of its target space is examined.  Moreover it
is found that any weak HKT manifold solves all the conditions required by
$N=4B$ one dimensional supersymmetry. The
understanding of HKT geometries requires the investigation of various
examples.  In this note we present a class of invariant
hypercomplex structures on nilpotent Lie groups which give rise  always to
weak HKT structures.  Moreover, they are Obata flat
when restricted to the 2-step nilpotent case (\cite{DF1}) but not in
general (see Example 3.3 in \cite{bilbao}).  In the
8-dimensional case there are only 3 possible groups and in dimension $4k,
k>2$ there are continuous families of weak HKT
structures (see Section 4). The HKT metrics in dimension 8 and their
properties are well understood (see 5.1). We note that
they have directions of positive Ricci curvature and directions of negative
Ricci curvature \cite{Mi} and their geodesics can
be given explicitly \cite{Ka, BTV}. It would be of interest to understand
the geodesic behaviour on compact quotients.  In
\cite{MS}  J. Michelson and A.Strominger proved that any weak HKT structure
which is quaternionic integrable (equivalently
Obata flat) can be constructed from a potential $L$ and ask whether
generically, one can do without the integrability
condition.  In particular all HKT structures on the 2-step case considered
in this note are associated to a potential.  It
would be of interest to see if the 12-dimensional example given in
\cite{bilbao} having a weak non integrable HKT structure can
be constructed from a potential.

\bigbreak 
\renewcommand{\thebibliography}{\list{\arabic{enumi}.\hfil}
{\settowidth\labelwidth{18pt}\leftmargin\labelwidth\advance
\leftmargin\labelsep\usecounter{enumi}}\def\newblock{\hskip.05em}
\sloppy\sfcode`\.=1000\relax}\newcommand{\bi}{\vspace{-3pt}\bibitem}

\end{document}